\documentclass[reqno]{amsart}

\theoremstyle{plain}
\newtheorem{theorem}{Theorem}[section]
\newtheorem{proposition}[theorem]{Proposition}
\newtheorem{lemma}[theorem]{Lemma}

\numberwithin{equation}{section}

\newcommand{\seclabel}[1]{\label{sec:#1}} 
\newcommand{\thmlabel}[1]{\label{thm:#1}} 
\newcommand{\lemlabel}[1]{\label{lem:#1}} 
\newcommand{\proplabel}[1]{\label{prop:#1}} 
\newcommand{\eqlabel}[1]{\label{eq:#1}} 

\newcommand{\secref}[1]{\ref{sec:#1}} 
\newcommand{\thmref}[1]{\ref{thm:#1}} 
\newcommand{\lemref}[1]{\ref{lem:#1}} 
\newcommand{\propref}[1]{\ref{prop:#1}} 
\renewcommand{\eqref}[1]{\ref{eq:#1}} 
\newcommand{\peqref}[1]{(\eqref{#1})} 

\newcommand\MM{\mathcal{M}}
\newcommand\QQ{\mathcal{Q}}
\newcommand\RR{\mathcal{R}}
\newcommand\LL{\mathcal{L}}
\newcommand\BB{\mathcal{B}}
\newcommand{\ld}{\backslash} 
\newcommand{\rd}{/}  

\title[Rectangular Quasigroups and Rectangular Loops]
{Rectangular Quasigroups and Rectangular Loops}

\author[M.~K.~Kinyon]{Michael~K.~Kinyon}
\address{Department of Mathematical Sciences \\
Indiana University South Bend \\
South Bend, IN 46634 USA}
\email{mkinyon@iusb.edu}
\urladdr{http://mypages.iusb.edu/\symbol{126}mkinyon}
\author[J.~D.~Phillips]{J.~D.~Phillips}
\address{Department of Mathematics \& Computer Science \\
Wabash College \\
Crawfordsville, IN 47933 USA}
\email{phillipj@wabash.edu}
\urladdr{http://www.wabash.edu/depart/math/faculty.html{\#}Phillips}

\date{\today}


\begin{document}

\begin{abstract}
We solve two problems posed by Krape\v{z} by finding a
basis of seven independent axioms for the variety of rectangular
loops. Six of these axioms form a basis for the variety of
rectangular quasigroups. The proofs of the lemmas showing
that the six axioms are sufficient are based on proofs
generated by the automated reasoning program OTTER,
while most of the models verifying the
independence of the axioms were generated by the finite
model builder Mace4.
\end{abstract}

\maketitle

\section{Introduction}
\seclabel{intro}

A \emph{rectangular band} is a semigroup $\BB = (B;\;\cdot )$, which is
is a direct product $\BB = \LL\oplus \RR$ of a left zero semigroup
$\LL = (L;\;\cdot )$ (that is, a semigroup satisfying $x\cdot y = x$) and
a right zero semigroup $\RR = (R;\;\cdot )$ (that is, a semigroup satisfying
$x\cdot y = y$). A semigroup that is a direct product of a group and
a rectangular band is called a \emph{rectangular group}. 

In a series of papers \cite{Kr1,Kr2,Kr3,Kr4}, Krape\v{z} has generalized
the notion of rectangular group to rectangular quasigroup and rectangular 
loop. Recall that a \emph{quasigroup} $\QQ = (Q;\;\cdot)$ is a set $Q$
together with a binary operation $\cdot : Q\times Q\to Q$ such that for
each $a\in Q$, the mappings $x \mapsto a\cdot x$ and
$x\mapsto x\cdot a$ are bijections. This induces two other
binary operations $\ld, \rd : Q\times Q\to Q$ as follows.
For $x,y\in Q$, $x\ld y$ is the unique solution $u$ to the
equation $x\cdot u = y$ and $x\rd y$ is the unique solution $v$
to the equation $v\cdot y = x$. This
leads to an equivalent definition: a \emph{quasigroup}
$\QQ = (Q;\;\cdot, \ld, \rd)$ is a set $Q$ with three binary operations
$\cdot, \ld, \rd : Q\times Q\to Q$ satisfying the equations:
\[
\begin{array}{c}
x\ld (x\cdot y) = y \qquad\qquad (x\cdot y) \rd y = x\\
x\cdot (x\ld y) = y  \qquad\qquad (x\rd y)\cdot y = x
\end{array}
\]
A \emph{loop} is usually defined as a quasigroup with a neutral
element $1\in Q$ satisfying $1\cdot x = x$ and $x\cdot 1 = x$. Loops
can also be characterized equationally as quasigroups satisfying the
additional axiom $x\ld x = y\rd y$. Basic references for quasigroup
and loop theory are \cite{Bel, Br, CPS, Pf}.

Following Krape\v{z}, a \emph{rectangular quasigroup} is a direct product
of a quasigroup and a rectangular band, while a \emph{rectangular loop} is
a direct product of a loop and a rectangular band. Both can be viewed
as varieties of algebras with three binary operations $\cdot$, $\ld$, $\rd$.
To consider a rectangular band as an algebra with three operations, one
simply sets $x\ld y = x\rd y = x\cdot y$. (This is a convention;
different choices lead to different axioms.)

In \cite{Kr1}, Krape\v{z} found a set of 12 axioms characterizing rectangular
loops as algebras $\MM = (M;\; \cdot, \ld, \rd)$. He indicated that he had shown
that some of the axioms are independent, but believed that not all of them are. He
also posed the problem of finding an independent set of axioms
(\cite{Kr1}, Problem 1, p. 66).
In \cite{Kr2}, he found a set of 15 axioms characterizing rectangular quasigroups,
and again posed the problem of finding an independent set of axioms. The two
axiom systems are quite distinct, that is, none of Krape\v{z}'s axioms for rectangular
quasigroups occur as axioms for rectangular loops.

In this paper, we solve both of the posed problems. We give a system of six axioms
and show that this is sufficient to characterize the variety of rectangular 
quasigroups. We then show that these same six plus one additional axiom characterize
the variety of rectangular loops. Finally, we will present models that show that
the seven rectangular loop axioms are independent. Here is the statement of our
main result.

\begin{theorem}
\thmlabel{main}
\begin{enumerate}
\item The equations
\[
\begin{array}{crclccrcl}
(Q1) & x\ld (xx) &=& x & &
(Q2) & (xx)\rd x &=& x \\
(Q3) & x (x\ld y) &=& x\ld (xy) & &
(Q4) & (x\rd y) y &=& (xy)\rd y \\
(Q5) & (x\ld y)\ld ((x\ld y)\cdot zu) &=& (x\ld (xz))u & &
(Q6) & (xy\cdot (z\rd u))\rd (z\rd u) &=& x((yu)\rd u) \\
\end{array}
\]
are an independent system of axioms for the
variety of rectangular quasigroups.

\item Equations (Q1)-(Q6) together with
\[
\begin{array}{crcl}
(L) & x\ld (x (y\ld y)) &=& ((x\rd x)y)\rd y 
\end{array}
\]
are an independent system of axioms for the
variety of rectangular loops.
\end{enumerate}
\end{theorem}

If part (1) of Theorem \thmref{main} is assumed, then it is easy to show
that (L), (Q1)-(Q6) characterize the variety of rectangular loops.
Indeed, if $\MM = (M;\; \cdot, \ld, \rd)$ satisfies (L), (Q1)-(Q6), 
then $\MM$ is a rectangular quasigroup, and hence a direct product
$\MM = \BB \oplus \QQ$ of a rectangular band $\BB$ and a 
quasigroup $\QQ$, each of which satisfies (L). This implies that
$\QQ$ is a loop, and so $\MM$ is a rectangular loop. Conversely,
since every rectangular band and every loop satisfy (L), so does every
rectangular loop. Thus what remains is to show that (Q1)-(Q6)
characterize the variety of rectangular quasigroups (\S\secref{proofs})
and then to show the independence of (L), (Q1)-(Q6) (\S\secref{indep}).

As part of the proof of Theorem \thmref{main}, we will also show
that eight of Krape\v{z}'s fifteen axioms are sufficient
to characterize rectangular quasigroups.

We also thank A. Krape\v{z} for helpful comments on the first draft of this paper.

\section{OTTER and Mace4}
\seclabel{otter-mace}

The proofs of the lemmas in \S\secref{proofs} are based on
proofs generated by the equational reasoning program OTTER
developed by McCune \cite{Otter}. OTTER can prove theorems
from axioms in first-order logic, but is strongest in equational
reasoning. Not surprisingly then, most \textit{new} mathematics
to come out of automated reasoning has been in fields
close to algebra. For general methods for applying
automated reasoning to problems in mathematics and other areas,
see the book by Wos and Pieper \cite{wosp}.
New results proved by OTTER in particular can be found
in the book by McCune and Padmanabhan \cite{mcpad}.

It is mathematically sound to use OTTER output directly
as the proof of a theorem, as is the practice in
\cite{mcpad}, for instance. Despite the complexity of
OTTER's search procedure,
the program can be made to output a \textit{proof object},
which can be independently verified by other software, such
as a \verb+lisp+ program. However, OTTER's proofs are often long
sequences of unintuitive equations, and it is useful to re-express
them in a form which a human reader can easily verify.
Some discussion of the procedure for ``humanizing'' proofs
occurs in \cite{hk}. In loop and quasigroup theory, this was
first applied by Kunen \cite{kuna, kunb, kunc, kund, kune},
and then later by the present authors and Kunen in various
combinations of coauthors
\cite{diass, arif, note, dcc, bol, tri, extra}.

The proofs in \S\secref{proofs} of the present paper are
somewhat closer to the OTTER proofs than in the 
aforementioned references. In quasigroup and loop theory,
one can rely on a great deal of existing machinery 
(e.g., autotopisms) to simplify proofs. Not so much
machinery is available for rectangular quasigroups and
rectangular loops. Nevertheless, the proofs herein are
still translated into a humanly verifiable form. 
In particular, we took care to ensure that each step
directly uses one of the axioms or a closely related
equation, or else uses one of the equations
in Proposition \propref{krapez} below.

The models in \S\secref{indep} that show the independence of 
axioms (Q1)-(Q6) were found using McCune's finite model builder
Mace4 \cite{Mace}. Mace4 can generate its output in a portable form,
which can then be used by other programs to independently verify the
claimed properties of the models. However, in this case all of the
models are quite small, so it was just as easy to verify the properties
``by hand''.

\section{Proofs}
\seclabel{proofs}

Since every quasigroup and every rectangular band satisfy
(Q1)-(Q6), to show that these equations characterize the variety
of rectangular quasigroups, it is enough to show that they
imply the axioms of Krape\v{z}:

\begin{proposition}[\cite{Kr2}]
\proplabel{krapez}
The following equations axiomatize the variety of
rectangular quasigroups: (Q1)--(Q4) and
\begin{align*}
(xy) \ld (xy\cdot z) &= x \ld (xz) \tag{K5} \\
(x\ld y) \ld ((x\ld y) z) &= x \ld (xz) \tag{K6}\\
(x\rd y) \ld ((x\rd y) z) &= x \ld (xz) \tag{K7}\\
x (y \ld (yz)) &= xz \tag{K8} \\
((xy) \rd y)z &= xz \tag{K9}\\
(x \cdot yz) \rd (yz) &= (xz)\rd z \tag{K10}\\
(x (y\ld z)) \rd (y\ld z) &= (xz)\rd z \tag{K11}\\
(x (y\rd z)) \rd (y\rd z) &= (xz)\rd z \tag{K12}\\
x \ld (x ((yz)\rd z)) &= ((x\ld (xy))z)\rd z \tag{K13}\\
(x \ld (xy))z &= x \ld (x\cdot yz) \tag{K14}\\
x ((yz)\rd z) &= (xy\cdot z)\rd z \tag{K15}
\end{align*}
\end{proposition}

\begin{lemma}
\lemlabel{lem1}
(Q1), (Q3), (K14) $\implies$
\begin{equation}
\eqlabel{tmp}
x\ld (x\cdot xz) = x(x\ld (xz)) = x\cdot x(x\ld z) = xz
\end{equation}
\end{lemma}

\begin{proof}
Replace $y$ with $x$ in (K14) and use (Q1) to get
$x\ld (x\cdot xz) = xz$. The other forms of 
\peqref{tmp} follow applying (Q3).
\end{proof}

\begin{lemma}
\begin{enumerate}
\item (Q1), (Q3), (K6), (K14) $\implies$ (K5), (K8)
\item (Q2), (Q4), (K12), (K15) $\implies$ (K10), (K9)
\end{enumerate}
\end{lemma}

\begin{proof}
For (1): By Lemma \lemref{lem1}, we may use \peqref{tmp}.
In (K6), replace $y$ with $x\cdot xy$ and use \peqref{tmp} to get (K5). Next
\[
\begin{array}{rcll}
x(y\ld (yz)) &=& x\cdot y(y\ld z) & \text{by\ (Q3)} \\
&=& x\cdot x(x\ld [y(y\ld z)]) & \text{by\ \peqref{tmp}} \\
&=& x\{ x(x\ld y)\cdot ([x(x\ld y)]\ld [y(y\ld z)])\} & \text{by\ (K5)} \\
&=& x\{ [x(x\ld y)]\ld [x(x\ld y)\cdot y(y\ld z)]\} & \text{by\ (Q3)}\\
&=& x\{ [x(x\ld y)]\ld [(x\ld (xy))\cdot y(y\ld z)]\} & \text{by\ (Q3)} \\
&=& x\{ [x(x\ld y)]\ld [x\ld (x(y\cdot y(y\ld z))]\} & \text{by\ (K14)} \\
&=& x\{ [x(x\ld y)]\ld [x\ld (x\cdot yz)]\} & \text{by\ \peqref{tmp}} \\
&=& x\{ [x(x\ld y)]\ld [(x\ld (xy))z)]\} & \text{by\ (K14)} \\
&=& x\{ [x(x\ld y)]\ld [x(x\ld y)\cdot z]\} & \text{by\ (Q3)}\\
&=& x\{ x(x\ld y)\cdot [(x(x\ld y))\ld z]\} & \text{by\ (Q3)} \\
&=& x\cdot x(x\ld z) & \text{by\ (K5)} \\
&=& xz & \text{by\ \peqref{tmp}}
\end{array}
\]
This is (K8).

The proof of (2) is the mirror of that of (1).
\end{proof}

\begin{lemma}
\begin{enumerate}
\item (Q4), (K5), (K9) $\implies$ (K7)
\item (Q3), (K10), (K8) $\implies$ (K11)
\end{enumerate}
\end{lemma}

\begin{proof}
For (1): 
\[
\begin{array}{rcll}
(x\rd y) \ld [(x\rd y)z] &=& [(x\rd y)y]\ld [(x\rd y)y\cdot z] & \text{by\ (K5)}\\
&=& [(x\rd y)y\cdot z]\ld [((x\rd y)y\cdot z)z] & \text{by\ (K5)} \\
&=& [((xy)\rd y)z]\ld [((xy)\rd y)z\cdot z] & \text{by\ (Q4)\ twice} \\
&=& (xz)\ld (xz\cdot z) & \text{by\ (K9)\ twice}\\
&=& x\ld (xz) & \text{by\ (K5)}
\end{array}
\]
which establishes (K7).

The proof of (2) is the mirror of that of (1).
\end{proof}

\begin{lemma}
(Q1), (Q3), (K14), (K15) $\implies$ (K13).
\end{lemma}

\begin{proof}
By Lemma \lemref{lem1}, we may use \peqref{tmp}. We compute
\[
\begin{array}{rcll}
[(x\ld (xy))z]\rd z &=& [x(x\ld y)\cdot z]\rd z & \text{by\ (Q3)} \\
&=& x[((x\ld y)z)\rd z] & \text{by\ (K15)} \\
&=& x\ld (x \{ x[((x\ld y)z)\rd z]\} ) & \text{by\ \peqref{tmp}} \\
&=& x\ld (x\{ [x(x\ld y)\cdot z]\rd z\} ) & \text{by\ (K15)} \\
&=& x\ld ([(x\cdot x(x\ld y)) z]\rd z) & \text{by\ (K15)} \\
&=& x\ld ([xy\cdot z]\rd z) & \text{by\ \peqref{tmp}} \\
&=& x\ld (x[(yz)\rd z]) & \text{by\ (K15)}
\end{array}
\]
and this is (K13).
\end{proof}

Combining Proposition \propref{krapez} and the lemmas, we obtain
the following.

\begin{theorem}
\thmlabel{krapez2}
Krape\v{z}'s axioms (Q1)-(Q4), (K6), (K12), (K14), (K15)
are sufficient to characterize the variety of rectangular
quasigroups.
\end{theorem}

So to complete the characterization part of the proof of
Theorem \thmref{main}, we need only the following.

\begin{lemma}
\begin{enumerate}
\item (Q1), (Q3), (Q5) $\implies$ (K6), (K14).
\item (Q2), (Q4), (Q6) $\implies$ (K12), (K15).
\end{enumerate}
\end{lemma}

\begin{proof}
For (1): First, we compute
\[
\begin{array}{rcll}
(x\ld y)\ld [(x\ld y)z] &=& (x\ld y)\ld [(x\ld y)(z\ld (zz))] & \text{by\ (Q1)} \\
&=& (x\ld y)\ld [(x\ld y)\cdot z(z\ld z)] & \text{by\ (Q3)} \\
&=& (x\ld (xz))(z\ld z) & \text{by\ (Q5)}
\end{array}
\]
so that the expression $(x\ld y)\ld [(x\ld y)z]$ is constant in $y$,
that is,
\[
(x\ld y)\ld [(x\ld y)z] = (x\ld u)\ld [(x\ld u)z] .
\]
Take $u = xx$ and apply (Q1) to get (K6). Next,
\[
\begin{array}{rcll}
x\ld (x\cdot yz) &=& (x\ld u)\ld [(x\ld u)\cdot yz] &\text{by\ (K6)}\\
&=& x\ld (xy\cdot z) & \text{by\ (Q5)}
\end{array}
\]
and this is (K14).

The proof of (2) is the mirror of that of (1).
\end{proof}


\section{Independence of the Axioms}
\seclabel{indep}

In this section, we present models that show the independence of 
axioms (L), (Q1)-(Q6), and this will complete the proof of the Theorem.
Note that equation (L) is equivalent to its own mirror, while the other axioms
come in mirrored pairs. Thus once we have presented a model satisfying,
for instance, all axioms except (Q6), it follows that the same underlying
set with the dual
operations $x \odot y := y\cdot x$, $x\ld\ld y := y\rd x$,
$x\rd\rd y := y \ld x$ will be a model satisfying all axioms except (Q5).
Thus it is enough to present four models.

The independence of (L) is obvious, because any nonloop quasigroup
satisfies (Q1)-(Q6), but not (L). Table 1 gives a specific example.
Table 2 is a model satisfying (L), (Q1), and (Q3)-(Q6), but not (Q2).
Table 3 is a model satisfying (L), (Q1)-(Q3), (Q5)-(Q6), but not (Q4).
Table 4 is a model satisfying (L), (Q1)-(Q5), but not (Q6).

\begin{table}[htb]
\label{table:notL}
\[
\begin{array}{lcccr}
\begin{array}{c|ccc}
\bullet &  0& 1& 2 \\
\hline
  0 &  0& 1& 2 \\
  1 &  2& 0& 1 \\
  2 &  1& 2& 0 \\
\end{array}
&&
\begin{array}{c|ccc}
\ld &  0& 1& 2 \\
\hline
  0 &  0& 1& 2 \\
  1 &  1& 2& 0 \\
  2 &  2& 0& 1 \\
\end{array}
&&
\begin{array}{c|ccc}
\rd &  0& 1& 2 \\
\hline
  0 &  0& 1& 2 \\
  1 &  2& 0& 1 \\
  2 &  1& 2& 0 \\
\end{array}
\end{array}
\]
\caption{(Q1)-(Q6), but not (L)}
\end{table}

\begin{table}[htb]
\label{table:notQ2}
\[
\begin{array}{lcccr}
\begin{array}{c|cc}
\bullet &  0& 1 \\
\hline
  0 &  0& 1 \\
  1 &  1& 0 \\
\end{array}
&&
\begin{array}{c|ccc}
\ld &  0& 1 \\
\hline
  0 &  0& 1 \\
  1 &  1& 0 \\
\end{array}
&&
\begin{array}{c|ccc}
\rd &  0& 1 \\
\hline
  0 &  1& 0 \\
  1 &  0& 1 \\
\end{array}
\end{array}
\]
\caption{(L), (Q1), (Q3)-(Q6), but not (Q2)}
\end{table}

\begin{table}[htb]
\label{table:notQ4}
\[
\begin{array}{lcccr}
\begin{array}{c|cccc}
\bullet &  0& 1& 2 & 3\\
\hline
  0 &  0& 1& 1& 0 \\
  1 &  0& 1& 1& 0 \\
  2 &  3& 2& 2& 3 \\
  3 &  3& 2& 2& 3 \\
\end{array}
&&
\begin{array}{c|cccc}
\ld &  0& 1& 2& 3 \\
\hline
  0 &  0& 1& 1& 0 \\
  1 &  0& 1& 1& 0 \\
  2 &  3& 2& 2& 3 \\
  3 &  3& 2& 2& 3 \\
\end{array}
&&
\begin{array}{c|cccc}
\rd &  0& 1& 2& 3 \\
\hline
  0 &  0& 2& 1& 0 \\
  1 &  0& 1& 1& 0 \\
  2 &  0& 2& 2& 0 \\
  3 &  3& 1& 1& 3 \\
\end{array}
\end{array}
\]
\caption{(L), (Q1)-(Q3), (Q5)-(Q6), but not (Q4)}
\end{table}

\begin{table}[htb]
\label{table:notQ6}
\[
\begin{array}{lcccr}
\begin{array}{c|cc}
\bullet &  0& 1 \\
\hline
  0 &  0& 1 \\
  1 &  0& 1 \\
\end{array}
&&
\begin{array}{c|cc}
\ld &  0& 1 \\
\hline
  0 &  0& 1 \\
  1 &  0& 1 \\
\end{array}
&&
\begin{array}{c|cc}
\rd &  0& 1 \\
\hline
  0 &  0& 0 \\
  1 &  0& 1 \\
\end{array}
\end{array}
\]
\caption{(L), (Q1)-(Q5), but not (Q6)}
\end{table}


\end{document}